\documentclass[a4,11pt]{amsart}%
\usepackage{amsmath,amssymb,amsfonts,enumerate,amsthm, amscd,}
\usepackage[all]{xy}
\usepackage{graphicx}
\usepackage{amsmath}
\usepackage{amsfonts}
\usepackage{amssymb}
\usepackage{hyperref}%
\usepackage{graphicx,enumitem,lipsum}
\usepackage{caption}
\usepackage{subcaption}
\usepackage{tikz-cd}
\usepackage{tikz}
\usetikzlibrary{calc}
\usetikzlibrary{decorations.pathmorphing}

\tikzset{curve/.style={settings={#1},to path={(\tikztostart)
			.. controls ($(\tikztostart)!\pv{pos}!(\tikztotarget)!\pv{height}!270:(\tikztotarget)$)
			and ($(\tikztostart)!1-\pv{pos}!(\tikztotarget)!\pv{height}!270:(\tikztotarget)$)
			.. (\tikztotarget)\tikztonodes}},
	settings/.code={\tikzset{quiver/.cd,#1}
		\def\pv##1{\pgfkeysvalueof{/tikz/quiver/##1}}},
	quiver/.cd,pos/.initial=0.35,height/.initial=0}

\tikzset{tail reversed/.code={\pgfsetarrowsstart{tikzcd to}}}
\tikzset{2tail/.code={\pgfsetarrowsstart{Implies[reversed]}}}
\tikzset{2tail reversed/.code={\pgfsetarrowsstart{Implies}}}

\tikzset{no body/.style={/tikz/dash pattern=on 0 off 1mm}}

\tikzset{no body/.style={/tikz/dash pattern=on 0 off 1mm}}
\providecommand{\U}[1]{\protect\rule{.1in}{.1in}}

\newcommand{\fm}{\mathfrak{m}}

\newtheorem{thm}{Theorem}[section]
\newtheorem{cor}[thm]{Corollary}
\newtheorem{lem}[thm]{Lemma}

\newtheorem{defn}[thm]{Definition}

\newtheorem{exam}[thm]{Example}

\theoremstyle{definition}

\theoremstyle{remark}

\theoremstyle{Definition and Notation}

\begin{document}

\title[On $\phi$-Pr\"ufer like conditions]{On $\phi$-Pr\"ufer like conditions }
\author{Adam Anebri}
\address{Adam Anebri\\ Euromed University of Fes, UEMF, Morocco
	$$E-mail\ address:\  a.anebri@ueuromed.org$$}

\author{Najib Mahdou}
\address{Najib Mahdou\\Department of Mathematics, Faculty of Science and Technology of Fez, Box 2202,
	University S.M. Ben Abdellah Fez, Morocco.
	$$E-mail\ address:\ mahdou@hotmail.com$$}

\author{El Houssaine Oubouhou}
\address{El Houssaine Oubouhou, Department of Mathematics, Faculty of Science and Technology of Fez, Box 2202,
	University S.M. Ben Abdellah Fez, Morocco.
	$$E-mail\ address:\ hossineoubouhou@gmail.com$$}

	
	\begin{abstract}  In this paper, we investigate  the question of when a $\phi$-ring is $\phi$-Pr\"ufer using two types of techniques: first, by analysing the lattice structure of the nonnil ideals of $\phi$-rings; and secondly,
by considering content ideal techniques which were developed to study Gaussian polynomials. In particular, we conclude that  every Gaussian $\phi$-ring is $\phi$-Pr\"ufer.  Key concepts such as  $\phi$-weak global dimension, primary ideals and  irreducible ideals  are discussed, along with their hereditary properties in $\phi$-Pr\"ufer rings.  We also prove that any semi-local  $\phi$-Pr\"ufer ring is a $\phi$-B\'ezout ring. This paper includes several theorems and examples that provide insights into the $\phi$-Pr\"ufer rings and their implications in the field of ring theory.
	\end{abstract}
	
	\maketitle

	
	\bigskip
	

	
	\section{Introduction}
	We assume throughout that all rings are commutative with  non-zero identity. Let $R$ be a ring, then $T(R)$ denotes the total quotient ring of $R$, $Nil(R)$ denotes the ideal of all nilpotent elements of $R$, and $Z(R)$ denotes the set of zerodivisors of $R$. We begin by recalling some background material. In 1932, Pr\"ufer introduced and conducted a study of integral domains in which any finitely generated ideal is invertible, as mentioned in \cite{prufer}. In 1936, Krull \cite{krull} named these rings after  Pr\"ufer and provided equivalent conditions that define a Pr\"ufer domain. Over the years, Pr\"ufer domains have been characterized by numerous equivalents, each of which can be extended in various ways to the context of rings with zerodivisors. In their paper dedicated to Gaussian properties, Bazzoni and Glaz demonstrated that a Pr\"ufer ring satisfies any of the other four Pr\"ufer conditions if and only if its total ring of quotients also satisfies the same condition, see \cite[Theorems 3.3, 3.6, 3.7 and 3.12]{BMM}.
		Recall that a non-zerodivisor element of a ring $R$ is called a regular element and an ideal of $R$ is said to be regular if it contains a regular element. A ring $R$ is called a Pr\"ufer ring, in the sense of \cite{griffin}, if every finitely generated regular ideal of $R$ is invertible, i.e., if $I$ is a finitely generated regular ideal of $R$ and $I^{-1}=\{x \in T(R) \mid x I \subset R\}$, then $I I^{-1}=R$. A Pr\"ufer domain is a Pr\"ufer ring and a homomorphic image of a Pr\"ufer domain is a Pr\"ufer ring. Many characterizations and properties of Pr\"ufer rings are stated in \cite{boisen, griffin, lucas}. For further study of Pr\"ufer domains and Pr\"ufer rings, we recommend \cite{ Huckaba, kap}. Recall from \cite{boisen} that a ring $R$
		is called a pre-Pr\"ufer ring if every proper homomorphic image of $R$ is a Pr\"ufer ring, i.e., if $R / I$ is a Pr\"ufer ring for each nonzero proper ideal $I$ of $R$. Note that  the class of Pr\"ufer rings and the class of pre-Pr\"ufer rings are not comparable under set inclusion  (cf. \cite{boisen}). A ring $R$ is called a B\'ezout ring, in the sense of \cite{Huckaba}, if every finitely generated regular ideal of $R$ is principal.\par

		Extending this concept to rings with zerodivisors leads to the creation of five classes of Pr\"ufer-like rings, with various homological aspects (\cite{BG} and \cite{G}). At this point, we make the following definition:
\begin{defn}

 Let $R$ be a ring.
 \begin{enumerate}
      \item[$(1)$]  $R$ is called semi-hereditary if every finitely generated ideal of $R$ is projective.
\item[$(2)$]  $R$ is said to have weak global dimension $\leq 1$ $($w.gl.$\operatorname{dim}(R) \leq 1)$ if every finitely generated ideal of $R$ is flat \cite{G, Gw}.
\item[$(3)$]  $R$ is said to be arithmetical if the lattice formed by its ideals is distributive \cite{fu}.
\item[$(4)$]  $R$ is called a Gaussian ring if for every $f, g \in R[X]$, one has the content ideal equation $c(f g)=c(f) c(g)$ \cite{Ts}.
 \end{enumerate}
\end{defn}

In the domain context, all the concepts described earlier coincide with the definition of a Pr\"ufer domain. However, in the more general context of arbitrary rings, these notions become distinct. In \cite{G}, Glaz presents examples that illustrate the differences between these concepts. The following diagram of implications summarizes the relations between them:
$$
\text { Semi-hereditary } \Rightarrow \text { weak global dimension } \leq 1 \Rightarrow \text { Arithmetical } \Rightarrow \text { Gaussian } \Rightarrow \text { Pr\"ufer }
$$
		
		An ideal $I$ of a ring $R$ is said to be a nonnil ideal if $I \nsubseteq N i l(R)$.	Recall from \cite{DB}  that a prime ideal $P$ of $R$ is called a divided
		prime ideal if  it is comparable to every ideal of $R$. Set $ \mathcal{H} = \{R \mid R$ is a commutative ring
		and $Nil(R)$ is a divided prime ideal of $R\}.$ If $R \in  \mathcal{H}$, then $R$ is called a $\phi$-ring. A $\phi$-ring is called a strongly $\phi$-ring if $Z(R)=Nil(R)$. Recall from \cite{DB} that for a ring $R \in \mathcal{H}$ with total quotient ring $T(R)$, the map $\phi: T(R) \rightarrow R_{N i l(R)}$ such that $\phi\left(\frac{a}{b}\right)=\frac{a}{b}$ for $a \in R$ and $b \notin Z(R)$ is a ring homomorphism from $T(R)$ into $R_{N i l(R)}$, and $\phi$ restricted to $R$ is also a ring homomorphism from $R$ into $R_{N i l(R)}$ given by $\phi(x)=\frac{x}{1}$ for each $x \in R$. It is clear that $\phi(R)$ is ring-isomorphic to the ring of quotient $R / \operatorname{ker}(\phi)$. Note that if $I$ is a nonnil ideal of a $\phi$-ring $R$, then $\phi(I)$ is a regular ideal of $\phi(R)$, and so
		a nonnil ideal $I$ of a $\phi$-ring $R$ is called $\phi$-invertible if $\phi(I)$ is an invertble ideal of $\phi(R)$.
		In 2004, Anderson and Badawi \cite{phip} extended the notion of Pr\"ufer domains to that of $\phi$-Pr\"ufer rings which are $\phi$-rings $R$ satisfying that each finitely generated nonnil ideal is $\phi$-invertible.
		Also, they characterized $\phi$-Pr\"ufer rings from the perspective of ring structures, which says that a $\phi$-ring $R$ is $\phi$-Pr\"ufer if and only if $R_{\mathfrak{m}}$ is a $\phi$-chained ring for any maximal ideal $\mathfrak{m}$ of $R$ if and only if $R / Nil(R)$ is a Pr\"ufer domain if and only if $\phi(R)$ is Pr\"ufer. In 2018,
		Zhao \cite{phif} gave a homological characterization of $\phi$-Pr\"ufer rings: a $\phi$-ring $R$ with $Z(R)=Nil(R)$
		is  $\phi$-Pr\"ufer if and only if each submodule of a $\phi$-flat module is $\phi$-flat, if and only if
		each nonnil ideal of $R$ is $\phi$-flat. Recently, in \cite{KMO}, Kim et al. showed that for a $\phi$-ring $R$, every ideal of $R$ is $\phi$-flat if and only $R$ is a
		$\phi$-Pr\"ufer ring with $Z(R) = Nil(R)$. Additional homological  characterizations of $\phi$-Pr\"ufer can be found in the interesting articles \cite{KMO, KMO3, phif, ZWZ}.\par
		In this paper, we focus on exploring the connections between different extensions of the $\phi$-Pr\"ufer ring concept. We establish that many properties and characterizations of $\phi$-Pr\"ufer rings are similar to those found in Pr\"ufer domains. Using the results of this study, we have the following implications, none of which are reversible:
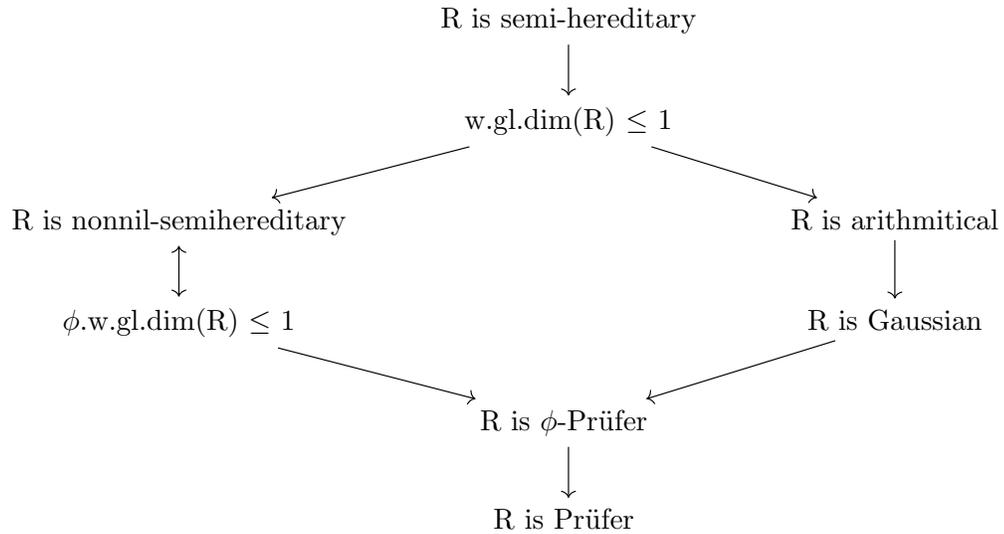
\[\begin{tikzcd}[ampersand replacement=\&]
	\& {\text{R is semi-hereditary}} \\
	\& {\text{w.gl.dim(R) $\leq$ 1}} \\
	{\text{R is nonnil-semihereditary}} \&\& {\text{R is arithmitical}} \\
	{\text{$\phi$.w.gl.dim(R) $\leq$ 1}} \&\& {\text{R is  Gaussian}} \\
	\& {\text{R is $\phi$-Pr\"ufer }} \\
	\& {\text{R is Pr\"ufer  }}
	\arrow[from=2-2, to=3-3]
	\arrow[from=5-2, to=6-2]
	\arrow[from=1-2, to=2-2]
	\arrow[from=3-3, to=4-3]
	\arrow[from=4-3, to=5-2]
	\arrow[from=4-1, to=5-2]
	\arrow[tail reversed, from=3-1, to=4-1]
	\arrow[from=2-2, to=3-1]
\end{tikzcd}\]
		In addition, when $R$ is a strongly $\phi$-ring, we discover several additional properties, and this leads to the construction of the following diagram:
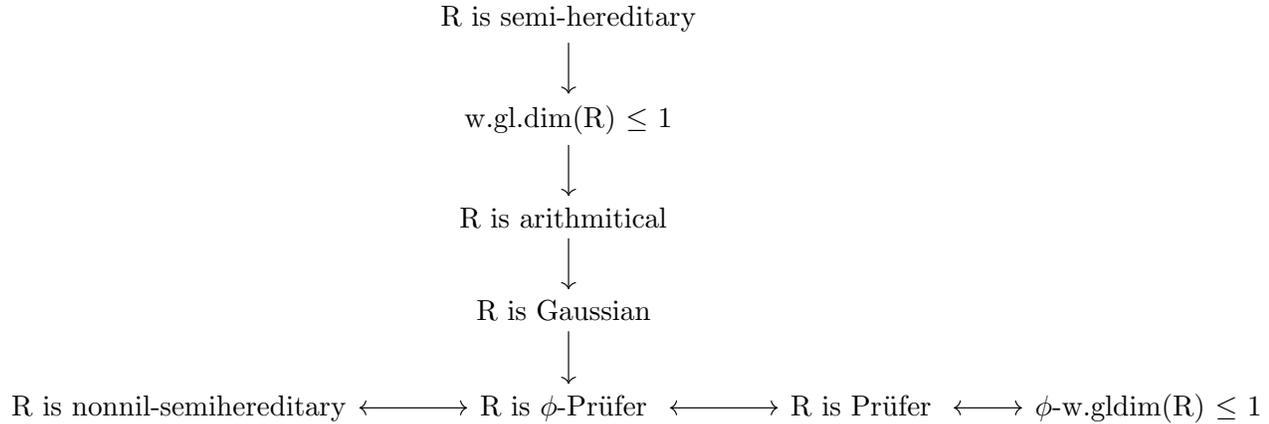
\[\begin{tikzcd}[ampersand replacement=\&]
	\& {\text{R is semi-hereditary}} \\
	\& {\text{w.gl.dim(R) $\leq$ 1}} \\
	\& {\text{R is arithmitical }} \\
	\& {\text{R is Gaussian }} \\
	{\text{R is nonnil-semihereditary}} \& {\text{R is $\phi$-Pr\"ufer  }} \& {\text{R is Pr\"ufer }} \& {\text{$\phi$-w.gldim(R) $\leq$ 1}}
	\arrow[from=2-2, to=3-2]
	\arrow[tail reversed, from=5-2, to=5-3]
	\arrow[from=3-2, to=4-2]
	\arrow[from=1-2, to=2-2]
	\arrow[tail reversed, from=5-2, to=5-1]
	\arrow[from=4-2, to=5-2]
	\arrow[tail reversed, from=5-3, to=5-4]
\end{tikzcd}\]
	Finally, we prove that in a $\phi$-Pr\"ufer ring $R$, every primary nonnil ideal of $R$ is irreducible and the converse is also true if any nonnil prime ideal of $R$ is maximal $($see Theorems \ref{pi}$)$. 
	\section{Main results}
Recall from \cite{phip} that a ring $R \in \mathcal{H}$ is called a $\phi$-chained ring if for every $x \in R_{N i l(R)} \backslash \phi(R)$, we have $x^{-1} \in \phi(R)$; equivalently, if for every $a, b \in R \backslash Nil(R)$, either $a \mid b$ or $b \mid a$ in $R$. We starts  by
the following lemmas.
\begin{lem}(\cite[ Theorem 2.7]{phip}) \label{l1}
 Let $R$ be a $\phi$-ring. Then $R$ is a $\phi-C R$ if and only if $R / N i l(R)$ is a valuation domain.
\end{lem}
\begin{lem} (\cite[Theorem 2.6]{phip}) \label{l2}
   Let $R$ be a $\phi$-ring. Then $R$ is a $\phi$-Pr\"ufer ring if and only if $R / N i l(R)$ is a Pr\"ufer domain.
\end{lem}

\begin{lem}(\cite[Theorem 2.9]{phip}) \label{l3} 
   Let $R$ be a $\phi$-ring. Then the following statements are equivalent:
\begin{enumerate}
    \item 
 $R$ is a $\phi$-Pr\"ufer ring;
   \item  $R_P$ is a $\phi-C R$ for each prime ideal $P$ of $R$;
   \item  $R_M$ is a $\phi-C R$ for each maximal ideal $M$ of $R$.
\end{enumerate}
\end{lem}

	\begin{thm} \label{t1}
		Let $R$ be a $\phi$-ring. Then the following assertions are equivalent:
		
		\begin{enumerate}
			\item[$(1)$]  $R$ is a $\phi$-Pr\"ufer ring. 
					\item[$(2)$] The lattice formed by the nonnil ideals of $R$ is distributive, i.e., 
$$
I \cap(J+K)=I \cap J+I \cap K
$$
for any three nonnil ideals of $R$.

		\end{enumerate}
	\end{thm}
	\proof
 $(1) \Rightarrow (2)$  Let $I$, $J$ and $K$ be three nonnil ideals of $R$. Since $R$ is a $\phi$-Pr\"ufer ring, $R/Nil(R)$ is a Pr\"ufer domain  by Lemma \ref{l2} and hence  The lattice formed by the ideals of $R/Nil(R)$ is distributive. Hence 
$$
{I/Nil(R)} \cap({J/Nil(R)}+{K/Nil(R)})={I/Nil(R)} \cap {J/Nil(R)}+{I/Nil(R)} \cap {K/Nil(R)}
$$
 and so $$
({I} \cap({J}+{K}))/Nil(R)=({I} \cap {J}+{I} \cap {K})/Nil(R)
.$$
We know $I, J$, and $K$ are not contained in $\operatorname{Nil}(R)$, and $\operatorname{Nil}(R)$ is prime, so both $I \cap(J+K)$ and $I \cap J+I \cap K$ cannot be contained in $\operatorname{Nil}(R)$. Since $\operatorname{Nil}(R)$ is divided, we know $\operatorname{Nil}(R)$ must be contained in both $I \cap(J+K)$ and $I \cap J+I \cap K$. By the above equality, we have

$$
I \cap(J+K)=I \cap J+I \cap K
$$

Therefore, the lattice formed by the nonnil ideals of $R$ must be distributive.\par

		$(2) \Rightarrow (1)$ Let ${\mathfrak{m}}$ be a maximal ideal of $R$.
 Since an ideal $J$ in $R_\fm$ is nonnil if and only  if $J=I_\fm$  for some nonnil ideal $I$ of $R$  and the formation of sums and intersections of ideals is preserved by extension from $R$ to $R_{\mathfrak{m}}$,   we can easily conclude that the lattice formed by  nonnil ideals of  is  $R_\fm$ distributive.   Here it is clearly enough to show that for two arbitrary  non-nilpotent  elements $a$ and $b$ in  $R_\fm$ will either $a \mid b$ or $b \mid a$. In fact, since the  nonnil ideals are assumed to form a distributive lattice, we have
$$
aR=aR \cap(b+(a-b))R=aR \cap bR +aR \cap(a-b)R
$$
so that $a$ may be written in the form $a=t+(a-b) c$, where $t$ is an element in $aR \cap bR$ and $b c$ an element in $aR$. Now, if $c$ is a unit, $b$ is a multiple of $b c$ and thus belongs to $aR$. If $c$ is not a unit, $(1-c)$ must be a unit, since the ring was supposed to be local. Therefore $a$ is a multiple of $a(1-c)=t-b c$, which is an element of $bR$. This means that we have either $a \mid b$ or $b \mid a$. Hence $R_\fm$ is  a  $\phi$-chained ring for any maximal ideal $\fm$ of $R$. Therefore $R$ is a $\phi$-Pr\"ufer ring according to Lemma \ref{l3}. \qed
\begin{thm}\label{t2}
     Let $R$ be a $\phi$-ring. A necessary and sufficient condition for $R$ to be $\phi$-Pr\"ufer ring is that for any pair of nonnil ideals $I$ and $J$ of $R$, such that
$I \subseteq J$ and $J$ is finitely generated, 
there exists a nonnil ideal $K$ for which $I=J \cdot K$.
\end{thm}

\proof To prove the ``if'' assertion let us consider an arbitrary pair of nonnil ideals $I, J$ such that $ I \subseteq J$ and $J$ is finitely generated, in a $\phi$-Pr\"ufer ring $R$. We assert that $(I: J)$ can be used as $K$.
In fact, to show that $J \cdot(I:J)=I$ it suffices to show that the local components agree. Taking into account that $J$ is finitely generated, then by \cite[Lemma 1]{jensen}  we have 
\begin{equation}\label{1}
    R_{\mathfrak{m}}(J \cdot(I: J))=R_{\mathfrak{m}}J \cdot R_{\mathfrak{m}}(I:J)=R_{\mathfrak{m}} J \cdot\left(R_{\mathfrak{m}}I: R_{\mathfrak{m}} J \right) .
\end{equation}

Since $J$ is finitely generated and the nonnil ideals of $R_{\mathrm{m}}$ are totally ordered, so $R_{\mathfrak{m}} J$ must be a principal  nonnil ideal containing $R_{\mathfrak{m}} I$, and consequently the right side of $(\ref{1})$ is equal to $R_{\mathfrak{m}} I$.  Now, to establish the sufficiency, using Lemma \ref{l3} we only need to   prove that the nonnil ideals in the local ring $R_{\mathfrak{m}}$, for any maximal ideal $\mathfrak{m}$, are totally ordered if $R$ satisfies the condition in the theorem.
It is clearly enough to prove that for any two  elements $\left[r_1 / s_1\right],\left[r_2 / s_2\right]$ in $R_{\mathfrak{m}}\setminus Nil(R_{\mathfrak{m}})$, $r_1, r_2, s_1, s_2 \in R$,  $r_1, r_2 \notin Nil(R)$  and  $s_1, s_2 \notin \mathfrak{m}$, at least one is a multiple of the other. In $R$, we consider the  nonnil ideals $I=r_1R$ and $J=r_1R+ r_2R$. Since $J$ is a finitely generated ideal containing $I$, the assumption about $R$ involves the existence of an ideal $K$ such that $I=JK$. This means that there exist two elements $x$ and $y$ in $R$ for which
$r_1=r_1 x+r_2 y$ with $r_2 x \in r_1R$ and $r_2 y \in r_1R$.

\begin{equation}\label{2}
    \left[r_1 / 1\right]=\left[r_1 / 1\right] \cdot[x / 1]+\left[r_2 / 1\right][y / 1]
\end{equation}
where
\begin{equation}\label{3}
     \left[r_2 / 1\right] \cdot[x / 1] \in \left[r_1 / 1\right]R_{\mathfrak{m}}.
\end{equation}

If $[x / 1]$ is a unit in $R_{\mathfrak{m}}$ we have by $(\ref{3})$
$$
\left[r_1 / s_1\right]\mid \left[r_1 / 1\right]\mid |\left[r_2 / 1\right] \cdot[x / 1] \mid \left[r_2 / 1\right] \mid\left[r_2 / s_2\right] \text {. }
$$
If $[x / 1]$ is a non-unit, $([1 / 1]-[x / 1])$ will be a unit in the local ring $R_{\mathrm{m}}$, and so we get by $(\ref{2})$
$$
\left[r_1 / 1\right] \cdot([1 / 1]-[x / 1])=\left[r_2 / 1\right] \cdot[y / 1]
$$
which implies
$$
\left[r_2 / s_2\right]\mid \left[r_2 / 1\right] \mid \left[r_1 / 1\right]  \cdot([1 / 1]-[x / 1]) \mid \left[r_1 / 1\right] \mid \left[r_1 / s_1\right] .
$$
Thus at least one of the elements $\left[r_1 / s_1\right]$ and $\left[r_2 / s_2\right]$ is a multiple of the other. \qed\\

Let $I$ be an ideal of $R$, and let $J$ be a nonempty subset of $R$. The residual of $I$ by $J$ is defined as $(I: J) = \{x \in R: xJ \subseteq I\}$. Furthermore, for $J = (a) \subseteq R$, we prefer the notation $(I: a)$ instead of $(I: (a))$.

\begin{thm}\label{t3}
    
Let $R$  be a $\phi$-ring. Then the following conditions are equivalent:
\begin{enumerate}
\item[$(1)$]   $R$ is a $\phi$-Pr\"ufer ring.
\item[$(2)$]  $(I+J): K=I: K+J: K$ for arbitrary  nonnil ideals $I$ and $J$, and any finitely generated ideal $K$.

\item[$(3)$]   $K:(I \cap J)=K: I+K: J$ for any finitely generated nonnil ideals $I$ and $J$ and arbitrary $K$.
\end{enumerate}
\end{thm} 
\proof 
$(1) \Rightarrow (2)$ Suppose that $R$ is $\phi$-Pr\"ufer. To prove the identity in (2) it suffices to show that the extensions to $R_{\mathfrak{m}}$ agree for any maximal ideal $\mathfrak{m}$. By \cite[Lemma 1]{jensen} and the assumption about $K$, sums and quotients of ideals are preserved by this extension, so that it will do to show the relation (2) for the ideals of $R_\fm$. But in these rings (2) surely is true, since the set of nonnil ideals is totally ordered according to Lemma \ref{l3}.\par
$ (1) \Rightarrow (3)$  Assume that  $R$ is $\phi$-Pr\"ufer.  Before proving the identity in (3), we notice that in any ring we have
$$
K:(I \cap J) \supseteq K: I+K: J
$$

so that we need only prove the converse inclusion. $I$ and $J$ are finitely generated nonnil ideals, but we do not know if $I \cap J$ is finitely generated, but anyway the \cite[Lemma 1]{jensen} implies
$$
R_{\mathfrak{m}}(K:(I \cap J)) \subseteq R_{\mathfrak{m}} K: R_{\mathfrak{m}}(I \cap J)=R_{\mathfrak{m}} K:\left(R_{\mathfrak{m}} I \cap R_{\mathfrak{m}} J\right)
$$
 and  so 
\begin{equation}\label{4}
     R_{\mathfrak{m}}(K: I+K: J)=R_{\mathfrak{m}} K: R_{\mathfrak{m}} I+R_{\mathfrak{m}} K: R_{\mathfrak{m}} J=R_{\mathfrak{m}} K:\left(R_{\mathfrak{m}} I \cap R_{\mathfrak{m}} J\right)
\end{equation}

the last equality in $(\ref{4})$ following from the fact that the set of nonnil ideals of $R_{\mathfrak{m}}$ is totally ordered by Lemma \ref{l3}. By the localization principle the desired inclusion is readily obtained.\\
 $ (2) \Rightarrow (1)$  Let $R$ be a ring for which (2) holds. We only need to  prove that for any two non-nilpotent elements $\left[r_1 / s_1\right]$ and $\left[r_2 / s_2\right]$ in $R_\fm,$ where  $ \fm$ is a maximal ideal, at least one is a multiple of the other according to Lemma \ref{l3}.  Since $\left[r_1 / s_1\right]$ and $\left[r_2 / s_2\right]$ are  non-nilpotent elements  of $R_{\mathfrak{m}}$, we get $r_1,r_2\in R\setminus Nil(R)$. If in (2) we choose $I=r_1R, J=r_2R, K=r_1R+ r_2R$, we get
$$
R=\left(r_1R: r_2R\right)+\left(r_2R: r_1R\right)
$$
so that there exist elements $x$ and $y$ for which
$$
1=x+y, \quad r_1\left|r_2 x, \quad r_2\right| r_1 y
$$
At least one of these elements does not belong to $\mathfrak{m}$, for instance $x \notin \mathfrak{m}$. Then $[x / 1]$ is a unit in $R_{\mathfrak{m}}$, and consequently
$$
\left[r_1 / s_1\right]\mid \left[r_1 / 1\right]\mid\left[r_2 x / 1\right] \mid \left[r_2 / 1\right]\mid \left[r_2 / s_2\right] .
$$
$(3) \Rightarrow (1)$ This may be proved similarly as $(2) \Rightarrow(1)$.\qed

\begin{thm}\label{t4}
    Let $R$ be a $\phi$-ring. Then the following assertions are equivalent:
    \begin{enumerate}
    \item[$(1)$] $R$ is a $\phi$-Pr\"ufer ring.
    \item[$(2)$] $({I} \cap {J}) \cdot {K}={I} {K} \cap {J}$ for any three nonnil ideals of $R$.
\end{enumerate}     
\end{thm}

\proof By passing to the local quotient rings $R_{\mathfrak{m}}$, it is enough to prove that a local  $\phi$-ring $R$ is $\phi$-Pr\"ufer if and only if $\left(2\right)$ holds for the nonnil ideals of $R$.
By Lemma \ref{l3} the "only if" part is obvious.
The "if" part is proved by showing that for any two non-nilpotent elements $a$ and $b$ will either $a \mid b$ or $b \mid a$ using Lemma \ref{l3}, provided that the assertion $\left(2\right)$ holds.
Putting ${I}=aR, {J}=bR$ and ${K}=aR+bR$, the condition $\left(2\right)$  implies
$$
abR \subseteq {K} \cap {J}=({I} \cap {J}) \cdot {K}=(aR \cap bR) \cdot(aR+bR)
$$
which in connection with the trivial converse inclusion involves
$$
abR=(aR \cap bR) \cdot(aR+bR) .
$$ Hence there exist two elements $x$ and $y \in aR\cap bR$ such that
$$
a b=x a+y b
$$
and thus $1=x / b+y / a$, $x / b$ and $y / a \in R$.
Since $R$ is local,  at least one of the elements $x / b$ and $y / a$ must be a unit; if $x / b$ is a unit we have $x \mid b$, which gives that  $a \mid x$ and $ x\mid b$, i.e., $a \mid b$.
Similarly, if $y / a$ is a unit we conclude that $b \mid a$. This completes the proof. \qed\\

Let $f\in R[X]$, the content of $f$ is the ideal of $R$ generated by the coefficients of $f$.
The content of $f$ is usually denoted by $c(f)$.  Note that if $R$ is a $\phi$-ring and  $f(X)=a_0+a_1 X+a_2 X^2+\cdots+a_n X^n, a_i \in R$, is a non-nilpotent element of  $R[X]$, then at least  one of the $a_i\in R \setminus Nil(R)$. Consider $\Tilde{f}=b_0+b_1X+b_2 X^2+\cdots+b_n X^n$ with $b_i=a_i$ if $a_i$ is non-nilpotent and $b_i=0$ if not. Since $Nil(R)$ is divided, we can easily see that  $c(f)=c(\Tilde{f})$. In addition,  if $g$ is a non-nilpotent element of  $R[X]$, then $c(fg)=c(\Tilde{f}g)=c(\Tilde{f}\Tilde{g})$ and consequently $c(fg)=c(f)c(g)$ if and only if $c(\Tilde{f}\Tilde{g})=c(\Tilde{f})c(\Tilde{g})$.

    \begin{thm}\label{t5}
Let $R$ be a $\phi$-ring. Then the  following assertions are equivalent:
\begin{enumerate}
    \item[$(1)$]  $R$ is a $\phi$-Pr\"ufer ring.
    \item[$(2)$]   Every finitely generated  nonnil ideal of $R$ is locally principal.
    \item[$(3)$]   Every non-nilpotent  element  $f$ in $R[X]$ is Gaussian.
    \item[$(4)$]   $c(fg)=c(f)c(g)$ for each non-nilpotent  elements  $f,g$ in $R[X]$
\end{enumerate}

\end{thm}
\proof
$(1)\Rightarrow (2)$ Assume that $R$ is a $\phi$-Pr\"ufer ring and let $I$ be a finitely generated nonnil ideal of $R$. So, we have $I_\fm$ is  nonnil finitely generated of $R_\fm$. On the other hand, since $R$ is a $\phi$-Pr\"ufer ring, we get $R_\fm$ is a $\phi$-chained ring by Lemma \ref{l3} and so $R_\fm/Nil(R_\fm)$ is a valuation domain. Hence $I_\fm +Nil(R_\fm)$ is a principal ideal of $R_\fm/Nil(R_\fm)$, whence $I_\fm$ is principal. 
Thus $I$ is locally principal.\par
  
$(2)\Rightarrow (3)$ Let  $f=a_0+ a_1X+\cdots +a_nX^n$ be a non-nilpotent element  of $R[X]$.  Let $\fm$ be a maximal ideal of $R$. As $c(f)$ is a nonnil ideal $R$, we have that $c(f)_\fm$ is a principal ideal of $R_\fm$   and it is generated by one of the $a_i$'s. Factorizing out that $a_i$ gives
$$
f=a_i \sum_j r_j X^j, \ \ \text{ with } r_j\in  R, r_i=1 .
$$
The element $\sum_j r_j  X^j$ has content $R$, hence is Gaussian by \cite[Lemma 2.7]{Ts}.
Let $g \in R[X]$. So, 
$$
\begin{aligned}
f g= & \left(a_i \sum r_j  X^j\right) g=a_i\left[\left(\sum r_j  X^j\right) g\right] \\
c(f g) & =\left(a_i\right) c\left(\left(\sum r_j  X^j\right) g\right) \\
& =\left(a_i\right) c(g) \text { since } \sum r_j  X^j \text { is Gaussian } \\
& =c(f) c(g) .
\end{aligned}
$$
 It follows that $f$ is Gaussian over every $R_{\fm}$   for every  maximal ideal  $\fm$ of $R$ and hence $f$ is Gaussian.\par
$(3)\Rightarrow (4)$ Straightforward.\par
 $(4) \Rightarrow (1)$ To prove that $R$ is a $\phi$-Pr\"ufer ring, it suffices to show that $R/Nil(R)$ is a  Pr\"ufer domain according to Lemma \ref{l2}.  Let $\Bar{f},\Bar{g} $ be two nonzero elements of  $ R/Nil(R)[X]=R[X]+Nil(R)[X]$, Hence $f$ and $g$ are non-nilpotent elements of $R[X]$. By assumption,  we obtain that $c(fg)=c(f)c(g)$, and by passing to the ring quotient and using the fact that $c(\Bar{f})=c(f)+Nil(R)[X]$, we can easily deduce that $c(\Bar{f}\Bar{g})=c(\Bar{f})c(\Bar{g})$. This yields that $R/Nil(R)$ is a Gaussian domain and consequently $R/Nil(R)$ is a Pr\"ufer domain.\qed\\

Combining the above five theorems with \cite[Corollary 2.10]{phip}, we arrive at the following corollary.
 
    \begin{cor}\label{cor0}
Let $R$ be a $\phi$-ring. Then the following statements are equivalent:
\begin{enumerate}
\item[$(1)$]  $R$ is a $\phi$-Pr\"ufer ring.
\item[$(2)$]  $\phi(R)$ is a Pr\"ufer ring.
\item[$(3)$]  $R / Nil(R)$ is a Pr\"ufer domain.
\item[$(4)$]  $\phi(R) / Nil(\phi(R))$ is a Pr\"ufer domain.
\item[$(5)$]  $R_P / Nil\left(R_P\right)$ is a valuation domain for each prime ideal $P$ of $R$.
\item[$(6)$]  $R_M / Nil\left(R_M\right)$ is a valuation domain for each maximal ideal $M$ of $R$.
\item[$(7)$]   The lattice formed by the nonnil ideals of $R$ is distributive, i.e., $$I \cap(J+K)=I \cap J+I \cap K  .$$  for any three nonnil ideals $I,J$ and $K$ of $R$.

\item[$(8)$]   For any pair of nonnil ideals $I$ and $J$ of $R$ such that
$I \subseteq J$ and $ J$ is finitely generated,  
there  exists a nonnil ideal $K$ satisfying $I=J \cdot K$.
\item[$(9)$]  $(I+J): K=I: K+J: K$ for arbitrary  nonnil ideals $I$ and $J$, and any finitely generated ideal $K$ of $R$.

\item[$(10)$]   $K:(I \cap J)=K: I+K: J$ for any finitely generated nonnil ideals $I$ and $J$ and an arbitrary  ideal $K$ of $R$.
\item[$(11)$]   $({I} \cap {J}) \cdot {K}={I} {K} \cap {J}$ for any three nonnil ideals of $R$.
 \item[$(12)$]   Every finitely generated  nonnil ideal of $R$ is locally principal.
  \item[$(13)$]   Every non-nilpotent  element  $f$ of $R[X]$ is Gaussian.
  \item[$(14)$]   $c(fg)=c(f)c(g)$ for any non-nilpotent  elements  $f,g$ of $R[X]$
\end{enumerate}
\end{cor}
\begin{cor}
Let $R$ be a $\phi$-ring. If $R$ is Gaussian, then $R$ is $\phi$-Pr\"ufer.
\end{cor}
The converse of the previous Corollary is not always true, as the following example indicates.

 \begin{exam}
     Let $K$ be a field and   let $R=K[x,y] /\left(x^2, y^2\right)$. Then  it is easy to see that  $R$ is a $\phi$-von Neumann regular ring $($because of $R/Nil(R)\cong K)$, and so $R$ is a $\phi$-Pr\"ufer ring by \cite[Corollary 4.9]{MO}. But $R$ is not Gaussian,   the polynomial $x Z+y \in R[Z]$ is not Gaussian since $(x Z+y)(xZ- y)=0$ while $(x, y)^2=(xy) \neq(0)$.
 \end{exam}

Let $R$ be a $\phi$-ring. Recall from \cite{strong} that the $\phi$-weak global dimension of $R$ is determined by the formulas:
$$
\begin{aligned}
\phi\text{-}w.g l . d i m(R) & =\sup \left\{f d_R(R / I) \mid I \text { is a nonnil ideal of } R\right\} \\
& =\sup \left\{f d_R(R / I) \mid I \text { is a finitely generated nonnil ideal of } R\right\} ,
\end{aligned}
$$
 and the  $\phi$-global dimension of $R$ is determined by the formula:
$$\phi\text{-}gl.dim(R) = \sup\{pd_R(R/I) \mid I \text{ is a nonnil ideal of } R\}.$$
Note that $\phi$-$w.g l . d i m(R)\leq \phi $-$gl.dim(R) $.
	We now turn our attention to rings of small $\phi$-weak dimension.  Recall from \cite{MO} that a $\phi$-ring $R$ is called $\phi$-von Neumann regular ring if every $R$-module is $\phi$-flat; equivalently $(R,Nil(R))$ is a local ring.
 The rings $R$ of $\phi$-$w.gl.\operatorname{dim} R=0$ are precisely the $\phi$-von Neumann regular rings \cite[Theorem 2.8]{strong} and as such $\phi$-Pr\"ufer. The rings $R$ of $\phi$-$w.gl.\operatorname{dim} R\leq 1$ are precisely the  $\phi$-Pr\"ufer strong $\phi$-ring \cite[Theorem 2.9]{strong}.  Recall that a ring $R$ is called nonnil-semihereditary if every finitely
generated nonnil ideal of $R$ is projective. See for instance \cite{KMO2,strong}.
 
 \begin{thm}
     Let $R$ be a $\phi$-ring. Then the following statements are equivalent:
     \begin{enumerate}
         \item[$(1)$]  $R$ is a $\phi$-Pr\"ufer strong $\phi$-ring.
         \item[$(2)$]   $\phi$-$w.gl.\operatorname{dim} R\leq 1$.
         \item[$(3)$]  $R$ is a nonnil-semihereditary ring.
     \end{enumerate}
 \end{thm}
 \proof
$(1) \Leftrightarrow(2)$ See \cite[Theorem 2.9]{strong}.\\
$(2)\Rightarrow (3)$ Assume that $\phi$-$w.gl.\operatorname{dim} R\leq 1$. Since  $\phi$-$w.gl.\operatorname{dim} R=\sup \{f d_R(R / I) \mid I$  is a finitely generated nonnil ideal of $ R\}\leq 1$, we get $f d_R(R / I)\leq 1$ for any nonnil ideal of $R$ and so every nonnil ideal of $R$ is flat. Hence every ideal of $R$ is $\phi$-flat by \cite[Proposition 2.6]{MO}. Consequently,  $R$ is a  strongly $\phi$-Pr\"ufer ring according to \cite[Corollary 3.2]{KMO}, and  hence $R/Nil(R)$ is a Pr\"ufer domain. In particular, $R/Nil(R)$ is a coherent  domain and $Z(R)=Nil(R)$, and whence $R$ is a nonnil-coherent ring by \cite[Corollary 3.2]{nc}. Since every finitely presented flat module is projective, we have that every finitely generated nonnil ideal $I$ of $R$ is projective.\\
$(3)\Rightarrow (2)$ It is obvious.\qed\\

 Our next results use the $A\ltimes M$ construction. Let $A$ be a ring and $M$ be an $A$-module. Then $A\ltimes M$, the {\it trivial} ({\it ring}) {\it extension of} $A$ {\it by} $M$, is the ring whose additive structure is that of the external direct sum $A \oplus M$ and whose multiplication is defined by $(r_1, m_1) (r_2, m_2) := (r_1r_2, r_1m_2+ r_2m_1)$ for all $r_1,r_2 \in A$ and all $m_1, m_2 \in M$. Mainly, trivial ring extensions have been useful for solving many open problems and conjectures in both commutative and non-commutative ring theory. See for instance \cite{BKM}.
\begin{exam}
    Let $R=\mathbb{Z}_2\ltimes \mathbb{Q}/\mathbb{Z}_2$. Then $R$ is a $\phi$-Pr\"ufer ring with  $\phi$-$w.gl.\operatorname{dim}(R)=\infty$.
   \end{exam} 
   \proof
    Consider $x=(2,0)$ and $y=(0,1/2)$. It is clear that $x \in Z(R)\setminus Nil(R)$,  $(0: x)=y R$ and $(0: y)=x R$.
Now, let $m_x$ (resp., $m_y$) denote the multiplication by $x$ (resp., $y$). Since $(0: x)=y R$ and $(0: y)=x R$, we have the following infinite flat resolution of $x R$ with syzygies $x R$ and $y R$
$$
\ldots \longrightarrow R \stackrel{m_y}{\longrightarrow} R \stackrel{m_x}{\longrightarrow} R \stackrel{m_y}{\longrightarrow} \ldots \stackrel{m_y}{\longrightarrow} R \stackrel{m_x}{\longrightarrow} R  \stackrel{\pi}{\longrightarrow} R/xR \longrightarrow 0.
$$
We claim that $x R$ and $y R$ are not flat. Indeed, recall that over a local ring, a projective module is necessarily free. So, no projective module is annihilated by $x$ or $y$. However, as $x R$ is annihilated by $y$ and $y R$ is annihilated by $x$, it follows that neither $xR$ nor $y R$ is projective. Further, $x R$ and $y R$ are finitely presented in view of the exact sequence $0 \rightarrow y R \rightarrow$ $R \rightarrow x R \rightarrow 0$. It follows that $x R$ and $y R$ are not flat (since a finitely presented flat module is projective). Hence $fd(R/xR)=\infty$ and consequently 
$\phi$-$w.gl.\operatorname{dim}(R)=\infty$.\qed\\

 Recall that a ring $R$ is called a B\'ezout ring if every finitely generated regular ideal of $R$ is principal, and a ring $R \in \mathcal{H}$ is a $\phi$-B\'ezout ring if $\phi(I)$ is a principal ideal of $\phi(R)$ for every finitely generated nonnil ideal $I$ of $R$; equivalently, if $R/Nil(R)$ is a B\'ezout domain see \cite{phip}.  
 In view of Theorem  \ref{t5}, it can be seen that any $\phi$-B\'ezout ring is a $\phi$-Pr\"ufer ring. The converse  is generally not true since a Pr\"ufer domain need not be a B\'ezout domain, however, we shall show that the converse will hold if some restriction is imposed on the ring. In fact, it turns out that in the semi-local case, i.e. provided the ring has only finitely many maximal ideals, we can prove the converse.
\begin{thm}\label{t11}  
    A semi-local $\phi$-Pr\"ufer ring is a $\phi$-B\'ezout ring. 
\end{thm}
\proof
 Let $\mathfrak{m}_1, \mathfrak{m}_2, \ldots, \mathfrak{m}_n$ be the finitely many maximal ideals of the  $\phi$-Pr\"ufer ring $R$. We have to show that any finitely generated nonnil ideal is principal, but it is clearly enough to show that any nonnil ideal $I=aR+bR$ generated by two elements of $R$ is principal.
Firstly, note that if $a$ (resp., $b$) is a nilpotent element of $R$, since $Nil(R)$ is divided and $bR$ (resp., $aR$) is a nonnil ideal of $R$, we get $a\in Nil(R)\subset bR$ (resp., $b\in Nil(R)\subset aR$ ), Thus $I=bR$ (resp., $I=aR$) is principal. Hence we can assume that $a$ and $b$ are non-nilpotent elements of $R$. So,
we will therefore construct  two elements $\alpha$ and $\beta$ in $R$ such that the local components of the principal ideal $(\alpha a+\beta b)R$ satisfy the conditions
 \begin{equation}\label{6}
     (\alpha a+\beta b) R_{\mathfrak{m}_i}=a R_{\mathfrak{m}_i}+b R_{\mathfrak{m}_i} \  \text{ for } 1 \leq i \leq n
 \end{equation} 
because the localization principle ensures that $(\alpha a+\beta b)=(a, b)$. Since the nonnil ideals of $R_{\mathfrak{m}_i}$ form a totally ordered set, for each $\mathfrak{m}_i$. Then  we have either $a R_{\mathfrak{m}_i} \subseteq b R_{\mathfrak{m}_i}$ or $b R_{\mathfrak{m}_i} \subseteq a R_{\mathfrak{m}_i}$.
Let us assume that the $\mathfrak{m}_i$'s are numbered such that $a R_{\mathfrak{m}_i} \subseteq b R_{\mathfrak{m}_i}$ for $1 \leq \mathrm{i} \leq k$ and $b R_{\mathfrak{m}_i} \subseteq a R_{\mathfrak{m}_i}$ for $k<i \leq n$. Since a maximal ideal is not contained in any other maximal ideal different from itself, hence  there exists for each $i$ an element $c_i$ such that $c_i \in \mathfrak{m}_i$, but $c_i \notin \mathfrak{m}_j$ for $i \neq j$. Set $\alpha=c_1 \cdots c_k$ and $\beta=c_{k+1} \cdots c_n$, then $\alpha \in \mathfrak{m}_i$ for $1 \leq i \leq k$, but $\alpha \notin \mathfrak{m}_i$ for $k<i \leq n$, and $\beta \in \mathfrak{m}_i$ for $k<i \leqq n$, but $\beta \notin \mathfrak{m}_i$ for $1 \leq i \leq k$. With this choice of $\alpha$ and $\beta$ we have for the principal ideal $(\alpha a+\beta b)R$
$$
\begin{gathered}
(\alpha a+\beta b) R_{\mathfrak{m}_i}=b R_{\mathfrak{m}_i} \text { for } \quad 1 \leq i \leq k, \\
(\alpha a+\beta b) R_{\mathfrak{m}_i}=a R_{\mathfrak{m}_i} \text { for } \quad k<i \leq n.
\end{gathered}
$$
Since in the first case $a R_{\mathfrak{m}_i} \subseteq b R_{\mathfrak{m}_i}$ and in the second $b R_{\mathfrak{m}_i} \subseteq a R_{\mathfrak{m}_i}$, in either case we have obtained $(\ref{6})$.\qed\\

It is a well-known fact that in a nonnil Noetherian ring, any irreducible nonnil ideal is primary \cite[Proposition 1.13]{hizem}, while a primary nonnil ideal need not be irreducible. It might be worth noticing that for a $\phi$-Pr\"ufer ring, the situation is just the opposite. In fact, in a valuation domain of rank $2$, all ideals are irreducible since they are totally ordered by set inclusion. However, not all of its ideals are primary. The statement that a primary ideal is irreducible will be provided in the following theorem.

	\begin{thm}\label{pi} Let $R$ be a $\phi$-ring.
	     If $R$ is a $\phi$-Pr\"ufer ring, then any primary nonnil ideal of $R$ is irreducible. The converse is true if any nonnil prime ideal of $R$ is maximal.
      \end{thm}	
\proof Let $\mathfrak{q}$ be a primary  nonnil ideal of $R$ with the prime ideal $\mathfrak{p}$ as its radical. Let us further assume that $\mathfrak{q}$ is represented as an intersection $q={I} \cap {J}$. We have to show that $\mathfrak{q}={I}$ or $\mathfrak{q}={J}$. By passage to the generalized quotient ring $R_{\mathfrak{p}}$, we get $\mathfrak{q} R_{\mathfrak{p}}={I} R_{\mathfrak{p}} \cap {J} R_{\mathfrak{p}}$. By Lemma \ref{l3}, the set of nonnil ideals of $R_{\mathfrak{p}}$ is totally ordered, so  $\mathfrak{q} R_{\mathfrak{p}}={I} R_{\mathfrak{p}}$ or $\mathfrak{q} R_{\mathfrak{p}}={J} R_{\mathfrak{p}}$. Suppose that $\mathfrak{q} R_{\mathfrak{p}}={I} R_{\mathfrak{p}}$. In this case, we shall finish the proof by showing that $\mathfrak{q}={I}$. Since $\mathfrak{q}$ is $\mathfrak{p}$-primary, the contraction of $\mathfrak{q} R_{\mathfrak{p}}={I} R_{\mathfrak{p}}$ to $R$ is $\mathfrak{q}$. The contraction of ${I} R_{\mathfrak{p}}$ to $R$ is the $S$-component ${I}_s$ of $I$, where $S$ denoting the complement of $\mathfrak{p}$ in $R$. Now, we have  clearly ${I} \subseteq {I}_S$ and  ${I}_S=\mathfrak{q}$, then it follows that  ${I} \subseteq \mathfrak{q}$. The converse of this inclusion is obvious. \\
Now, we will prove the converse under the additional hypothesis that every nonnil prime ideal is maximal.
Firstly, note that if $Nil(R)$ is  maximal then $R$ is  naturally a $\phi$-Pr\"ufer ring  and $R$ has no  nonnil proper ideal. In this case, the proof is clear. Now, we may assume that the set of nonnil prime ideals of $R$ is not empty. To obtain the "if" part, we have to show that the set of nonnil ideals in any $R_{\mathfrak{m}}, \mathfrak{m}$ maximal, is totally ordered by set inclusion. Now, $\mathfrak{m} R_{\mathfrak{m}}$ is the only nonnil prime ideal of $R_{\mathfrak{m}}$ and is therefore the radical of any  proper nonnil ideal in $R_{\mathfrak{m}}$. All nonnil ideals of $R_{\mathfrak{m}}$ are $\mathfrak{m} R_{\mathfrak{m}}$ primary, and consequently in a $1-1$ correspondence with the $\fm$-primary ideals of $R$. The intersection of two $\fm$-primary ideals is itself an $\fm$-primary ideal, so that the irreducibility of the $\fm$-primary ideals gives that these are totally ordered by set inclusion. Since the above $1-1$ correspondence is order-preserving, the set of nonnil ideals in $R_{\mathfrak{m}}$ is totally ordered, as desired. \qed


\begin{thebibliography}{99}

			\bibitem{phip} D. F. Anderson and A. Badawi, \textit{On $\phi$-Pr\"ufer rings and $\phi$-B\'ezout rings}, Houston J. Math. 30(2) (2004), 331--343.	
		\bibitem{nc} K. Bacem and B. Ali, \textit{Nonnil-coherent rings}, Beitr. Algebra Geom. 57(2) (2016), 297--305.
			\bibitem{DB} A. Badawi, \textit{On divided commutative rings}, Comm. Algebra 27(3) (1999), 1465--1474.
			
			\bibitem{nn} A. Badawi, \textit{On nonnil-Noetherian rings}, Comm. Algebra 31(4) (2003), 1669--1677.
			
			\bibitem{BKM} C. Bakkari, S. Kabbaj and N. Mahdou, \textit{Trivial extensions defined by Pr\"ufer conditions}, J. Pure Appl. Algebra 214(1) (2010), 53--60.
			\bibitem{BMM} C. Bakkari, N. Mahdou and H. Mouanis, \textit{Pr\"ufer-like conditions in subrings retract and applications}, Comm. Algebra 37(1) (2009), 47--55.
			\bibitem{BG} S. Bazzoni, S. Glaz, \textit{Pr\"ufer rings}, in: Multiplicative Ideal Theory in Commutative Algebra, Springer, (2006), 263--277.
			
			\bibitem{boisen} M. B. Boisen, J. R. and P. B. Sheldon, \textit{Pre-Pr\"ufer rings}, Pacific J. Math. 58 (1975), 331--344.
			\bibitem{BS}	H. R. Butts, W. Smith, \textit{Pr\"ufer rings}, Math. Z. 95 (1967), 196--211.
\bibitem{fu} L. Fuchs, \textit{Uber die Ideale arithmetischer Ringe}, Comment. Math. Helv. 23 (1949) 334--341.
			
\bibitem{G} S. Glaz, \textit{ Pr\"ufer conditions in rings with zero-divisors}, CRC Press Series of Lectures in Pure Appl. Math. 241 (2005), 272--282.
\bibitem{Gw} S. Glaz, \textit{The weak global dimension of Gaussian rings}, Proc. Amer. Math. Soc. 133(9) (2005), 2507--2513.
\bibitem{griffin} M. Griffin, \textit{Pr\"ufer rings with zerodivisors}, J. Reine Angew. Math. 240(1970), 55--67.

\bibitem{hizem} S. Hizem and A. Benhissi  \textit{Nonnil-Noetherian rings and the SFT property}, Rocky Mountain J. Math. 41 (2011), 1483--1500.
\bibitem{Huckaba} J. A. Huckaba, \textit{Commutative Rings with Zero Divisors}, Marcel Dekker, New York Basel, (1988).
\bibitem{jensen} G. U. Jensen, \textit{Arithmetical rings}, Acta Sci. Acad. Hungar. 17 (1966), 115--123.
\bibitem{kap} I. Kaplansky, \textit{Commutative Rings}, The University of Chicago Press, Chicago, (1974).
\bibitem{KMO}  H. Kim, N. Mahdou and E. H. Oubouhou, \textit{When every ideal is $\phi$-$P$-flat}, Hacet. J. Math. Stat. 52 (3) (2023), 708--720.
\bibitem{KMO3}  H. Kim, N. Mahdou and E. H. Oubouhou, \textit{Generalizations of Prüfer rings and Bézout rings}. São Paulo J. Math. Sci., 18 (2024), 126--141.
\bibitem{KMO2}  H. Kim, N. Mahdou and E. H. Oubouhou, \textit{On the $\phi$-weak global dimensions of polynomial rings and $\phi$-Pr\"ufer rings}, J. Algebra Appl., to appear.
\bibitem{krull} W. Krull,  \textit{Beitr\"age zur arithmetik kommutativer integrit\"atsebereiche}, Math. Z. 41 (1936), 545--577.
\bibitem{lucas} T. G. Lucas, \textit{Some results on Pr\"ufer rings}, Pacific J. Math. 124 (1986), 333--343.
\bibitem{MO}    N. Mahdou and E. H. Oubouhou,  \textit{On $\phi$-$P$-flat modules and $\phi$-von neumann regular rings}, J. Algebra Appl.,   23(09) (2024), 2450143.
			\bibitem{prufer} H. Pr\"ufer,  \textit{Untersuchungen \"uber Teilbarkeitseigenschaften in K\"orpern}, J. Reine Angew. Math. 168 (1932), 1--36.
	\bibitem{Ts} H. Tsang, \textit{Gauss’s Lemma}, Ph.D. Thesis, University of Chicago, Chicago, (1965).
			
	\bibitem{strong} X. L. Zhang \textit{Strongly $\phi $-flat modules, strongly nonnil-injective modules and their homology dimensions}, (2022), arXiv preprint arXiv:2211.14681.

\bibitem{phif} W. Zhao, \textit{On $\phi$-flat modules and $\phi$-Pr\"ufer rings},  J. Korean Math. Soc. 55(5) (2018), 1221--1233.		
		
	\bibitem{ZWZ} W. Zhao, F. Wang, and X. Zhang, \textit{On $\phi$-projective modules and $\phi$-Pr\"ufer rings},
			Comm. Algebra 48(7) (2020), 3079--3090.
			
		
	\end{thebibliography}
\end{document}